\documentclass[a4paper,11pt,twoside]{article}
\usepackage{amssymb,amsmath,latexsym,geometry,fancyhdr,lineno,hyperref,titletoc}
\contentsmargin{0pt}

\dottedcontents{section}[2em]{\vspace{-1mm}\small}{2em}{0pt}
\dottedcontents{subsection}[5em]{\vspace{-1mm}\small}{3em}{5pt}

\hypersetup{colorlinks,linkcolor= blue,citecolor=blue}
\newtheorem{theorem}{Theorem}[section]

\newtheorem{remark}[theorem]{Remark}
\newtheorem{claim}[theorem]{Claim}
\newtheorem{corollary}[theorem]{Corollary}
\newtheorem{definition}[theorem]{Definition}

\numberwithin{equation}{section}
\allowdisplaybreaks[0]

\begin{document}
\title{{\bf\Large A Liouville-type theorem for the coupled Schr\"odinger systems and the uniqueness of the sign-changing radial solutions}}
\author{\\
{ \textbf{\normalsize Haoyu Li}}\\%\footnote{}\\%\footnote{Corresponding author}
{\it\small Departamento de Matem\'atica,}\\
{\it\small Universidade Federal de S\~{a}o Carlos,}\\
{\it\small S\~{a}o Carlos-SP, 13565-905, Brazil}\\
{\it\small e-mail: hyli1994@hotmail.com}\\
\\
{ \textbf{\normalsize Ol\'impio Hiroshi Miyagaki}}\\%\footnote{}\\%\footnote{Corresponding author}
{\it\small Departamento de Matem\'atica,}\\
{\it\small Universidade Federal de S\~{a}o Carlos,}\\
{\it\small S\~{a}o Carlos-SP, 13565-905, Brazil}\\
{\it\small e-mail: olimpio@ufscar.br}}
\date{}
\maketitle
{\bf\normalsize Abstract.} {\small
In this paper, we study the sign-changing radial solutions of the following coupled Schr\"odinger system
\begin{equation}
    \left\{
   \begin{array}{lr}
     -{\Delta}u_j+\lambda_j u_j=\mu_j u_j^3+\sum_{i\neq j}\beta_{ij} u_i^2 u_j \,\,\,\,\,\,\,\,  \mbox{in }B_1 ,\nonumber\\
     u_j\in H_{0,r}^1(B_1)\mbox{ for }j=1,\cdots,N.\nonumber
   \end{array}
   \right.
\end{equation}
Here, $\lambda_j,\,\mu_j>0$ and $\beta_{ij}=\beta_{ji}$ are constants for $i,j=1,\cdots,N$ and $i\neq j$. $B_1$ denotes the unit ball in the Euclidean space $\mathbb{R}^3$ centred at the origin. For any $P_1,\cdots,P_N\in\mathbb{N}$, we prove the uniqueness of the radial solution $(u_1,\cdots,u_j)$ with $u_j$ changes its sign exactly $P_j$ times for any $j=1,\cdots,N$ in the following case: $\lambda_j\geq1$ and $|\beta_{ij}|$ are small for $i,j=1,\cdots,N$ and $i\neq j$.
New Liouville-type theorems and boundedness results are established for this purpose.
}

\medskip
{\bf\normalsize 2020 MSC:} {\small 35A02, 35B53, 35J47}

\medskip
{\bf\normalsize Key words:} {\small Coupled Schr\"odinger systems; uniqueness; sign-changing radial solutions; Liouville-type theorem.
}

\pagestyle{fancy}
\fancyhead{} % clear all header fields
\fancyfoot{} % clear all footer fields
\renewcommand{\headrulewidth}{0pt}
\renewcommand{\footrulewidth}{0pt}
\fancyhead[CE]{ \textsc{Haoyu Li \& Ol\'impio Hiroshi Miyagaki}}
\fancyhead[CO]{ \textsc{A Liouville-type theorem for the coupled Schr\"odinger systems and ...}}
\fancyfoot[C]{\thepage}
%%%%----------------------------------------------------------------------------------------------------------------------------
%%%%----------------------------------------------------------------------------------------------------------------------------

%%\tableofcontents
%%\pagewiselinenumbers
%%%%----------------------------------------------------------------------------------------------------------------------------
%%%%----------------------------------------------------------------------------------------------------------------------------

\section{Introduction}
\subsection{Main theorems}
This paper proves the uniqueness of the radial solutions that have a component-wise prescribed number of nodes for the following problem:
\begin{equation}\label{e:AAAA}
    \left\{
   \begin{array}{lr}
     -{\Delta}u_j+\lambda_j u_j=\mu_j u_j^3+\sum_{i\neq j}\beta_{ij} u_i^2 u_j \,\,\,\,\,\,\,\,  \mbox{in }B_1 ,\\
     u_j\in H_{0,r}^1(B_1)\mbox{ for }j=1,\cdots,N.
   \end{array}
   \right.
\end{equation}
Here, $\lambda_j,\mu_j>0$, and $\beta_{ij}=\beta_{ji}$ for $i,j=1,\cdots,N$ and $i\neq j$. $B_1$ is the unit ball in $\mathbb{R}^3$ centred at the origin. $H_{0,r}^1(B_1)$ is the subspace consisting of radial functions of $H_0^1(B_1)$.
To be precise, for any $P_1,\cdots,P_N\in\mathbb{N}$, we prove the uniqueness of the radial solution $(u_1,\cdots,u_N)$ to Problem (\ref{e:AAAA}) with $u_j$ changes its sign exactly $P_j$ times for any $j=1,\cdots,N$.
As proved in \cite{LiuWang2019,LiWang2021}, such a solution with prescribed number of nodes exists as long as $\beta_{ij}<B$ for any $i\neq j$ and for some $B=B(P_1,\cdots,P_N)>0$. Furthermore, Li and Wang \cite{LiWang2021} proved a infinitely many existence result for the above solutions.
For instance, in the case of the two-component systems, the theorem reads as
~~
\begin{itemize}
  \item []
  \emph{For Problem (\ref{e:AAAA}) with $N=2$, assume $\lambda_1=\lambda_2>0$ and $-\beta>\mu_1=\mu_2>0$. For any $P\in\mathbb{N}$, there are infinitely many radial solutions $(u_{k,1},u_{k,2})$ with $u_{k,1}$ and $u_{k,2}$ change their sign exactly $P$ times for any $k=1,2,\cdots$.}
\end{itemize}
~~
It is a natural to inquire about the boundedness or even the uniqueness of solutions \emph{within the regime of existence but outside of multiplicity}.
In this note,
we give a partial affirmative answer.
\begin{theorem}\label{t:main}
If $\lambda_1,\cdots,\lambda_N\geq1$, for any $P_1,\cdots,P_N\in\mathbb{N}$, there exists a positive number $b=b(\lambda_1,\cdots,\lambda_N;\mu_1,\cdots,\mu_N;P_1,\cdots,P_N)>0$ such that if $|\beta_{ij}|<b$, Problem (\ref{e:AAAA}) admits an unique solution $(u_1,\cdots,u_N)$ with $u_j(0)>0$ and $u_j$ changing its sign exactly $P_j$ times for $j=1,\cdots,N$.
\end{theorem}
The existence part has been proved in \cite{LiWang2021}. We only need to prove the uniqueness part.
We will apply the ideas in \cite{Ikoma2009} to the sign-changing radial solutions. A by-product is the following remark.
\begin{remark}\label{r:nondegeneracy}
Under the assumption of Theorem \ref{t:main}, the solution $(u_1,\cdots,u_N)$ is non-degenerate in $(H_{0,r}^1(B_1))^N$, i.e., for any $(\varphi_1,\cdots,\varphi_N)\in H_{0,r}^1(B_1)^N$, we have
\begin{align}
-\Delta \varphi_j +\lambda_j\varphi_j = 3\mu_j u_j^2\varphi_j + \sum_{i\neq j}\beta_{ij} u_i^2 \varphi_j+ 2\sum_{i\neq j}\beta_{ij} u_i u_j \varphi_i\nonumber
\end{align}
for $j=1,\cdots,N$.
\end{remark}

In order to prove the boundedness, the following Liouville-type theorem plays an important role. We are delighted to note that we benefit a lot from Professor Pavol Quittner, who provides the idea and various comments of the proof of Theorem \ref{t:Liouville}, cf. \cite{QuittnerCommunication}. Consider the following elliptic system.
\begin{equation}\label{e:CCCC}
    \left\{
   \begin{array}{lr}
     -{\Delta}u_j=\mu_j u_j^3+\sum_{i\neq j}\beta_{ij} u_i^2 u_j \,\,\,\,\,\,\,\,  \mbox{in }\mathbb{R}^3 ,\\
     u_j\in C^{2}(\mathbb{R}^3),\,u_j\mbox{ is radial},\\
     u_j\mbox{is zero or changes its sign at most }P_j\mbox{ times and }(u_1,\cdots,u_N)\neq0,\\
     |u_j|_{L^\infty(\mathbb{R}^3)}\leq1\mbox{ for }j=1,\cdots,N.
   \end{array}
   \right.
\end{equation}
Here, $\mu_j>0$ are constants for $j=1,\cdots,N$.
The following Liouville-type theorem holds.

\begin{theorem}\label{t:Liouville}
For any $P_1\cdots,P_N\in\mathbb{N}$, and any $\beta\geq0$, there exists a positive number $\varepsilon_0=\varepsilon_0(\mu_1,\cdots,\mu_N;P_1,\cdots,P_N;\beta)>0$ such that if $\beta_{ij}\in[-\varepsilon_0,\beta]$ for any $i,j=1,\cdots,N$ and $i\neq j$,
Problem (\ref{e:CCCC}) admits no non-trivial solution.
\end{theorem}
Theorem \ref{t:Lioville2} concerns the one spatial dimension case, which can be proved in a similar way. We will state it in Section 2.2. As an immediate application of Theorem \ref{t:Liouville} and Theorem \ref{t:Lioville2}, the following boundedness result is implied.
\begin{theorem}\label{t:boundedness}
For any $P_1\cdots,P_N\in\mathbb{N}$, and any $\beta\geq0$, there exist positive numbers $\varepsilon_0=\varepsilon_0(\lambda_1,\cdots,\lambda_N;\mu_1,\cdots,\mu_N;P_1,\cdots,P_N;\beta)>0$ and $C=C(\lambda_1,\cdots,\lambda_N;\mu_1,\cdots,\mu_N;P_1,\cdots,P_N;\beta)>0$ such that if $\beta_{ij}\in[-\varepsilon_0,\beta]$ for any $i,j=1,\cdots,N$ and $i\neq j$, the solution $U$ to
Problem (\ref{e:AAAA}) satisfies $|U|_\infty\leq C$.
\end{theorem}

\subsection{Some historical remarks}
In this note, we study the uniqueness of the sign-changing radial solutions to Problem (\ref{e:AAAA}) with component-wisely prescribed number of nodes. Regarding the uniqueness problem, the current studies mostly focus on the positive solutions, cf. \cite{Ikoma2009,YaoWei2011,LiuLiuChang2015,ZhouWang2020} and the references therein.
The studies on the uniqueness of the positive solution can be traced back to the analogue results for the scalar field counterpart of problem (\ref{e:AAAA}), i.e. the equation $-\Delta u+u=|u|^{p-2}u$. This problem has been studied extensively. See \cite{GidasNiNirenberg1979,Kwong1989,KabeyaTanaka1999}.
For the uniqueness of radial sign-changing solutions, the complete answer remains open. However, in the last decade some partial results have been proved. We refer to \cite{AoWeiYao2015} for a slightly subcritical problem and \cite{Tanaka2016} for several special cases.
In \cite{KormanLiSchmidt2012} the authors carried out their study with the help of computers.

In this note we give a uniqueness result (cf. Theorem \ref{t:main}) for the radial sign-changing solution to the coupled Schr\"odinger systems. In \cite{LiWang2021} it was proved that in certain regions of the parameters the solutions with prescribed number of nodes can be multiple. This is based on a combination of the parabolic equation and the Lusternik-Schnirelmann argument. However, there are still some cases where we can only prove existence but not multiplicity.
The aim of this note is to find some uniqueness results in the region of existence but outside the region of multiplicity. To this end, we need new Liouville-type theorems for the coupled elliptic systems. See Theorem \ref{t:Liouville} and Theorem \ref{t:Lioville2}. Since the seminal work \cite{GidasSpruck1981}, Liouville type theorems have been widely applied to the semilinear elliptic equations. There are also many results on the Liouville-type theorem for the elliptic systems. We quote \cite{DancerWeiWeth2010,QuitterSouplet2012,QuitternBook2019} and the references therein for the positive solutions and \cite{Quittner2021} for the sign-changing radial solutions in a very general case. It was shown by Quittner \cite{Quittner2021} that in addition to the number of nodes for each component, the number of nodes of the difference and sum of the components, i.e., the comparisons, are also crucial. In Theorem \ref{t:Liouville} and Theorem \ref{t:Lioville2} of this paper, we can find a special case where the information about the comparison can be omitted. However, as it is shown in Remark \ref{r:novelty}, in the wider regions, the comparisons are still essential.

\subsection{The idea of the prove and the organization of this note}
Our approach is based on an implicit function theorem argument analogue to \cite{Ikoma2009}. However, in his work \cite{Ikoma2009}, Ikoma focused on the positive solutions. In order to solve the problem concerning the sign-changing radial solutions, we should prove the corresponding a priori and the non-degeneracy of the sign-changing radial solutions. Therefore, our note is divided into the following parts:
\begin{itemize}
  \item In Section 2, we prove the corresponding Liouville-type theorems and the non-degeneracy of the radial sign-changing solution to the scalar field counterpart of Problem (\ref{e:AAAA});
  \item In Section 3, we prove Theorem \ref{t:boundedness} and Theorem \ref{t:main}.
\end{itemize}

\section{Preliminaries}
\subsection{Notations}
In this part, we introduce some notations, in particular the number of nodes of a continuous radial function $u:\Omega\to\mathbb{R}$. Here, we assume that $\Omega$ is radial, i.e., it is either a ball, an annulus, the exterior domain of a ball or the whole space $\mathbb{R}^N$.
\begin{definition}
For a continuous radial function $u:\Omega\to\mathbb{R}$, we define the number of nodes of the function $u$ to be the the largest number $k$ such that there exists a sequence of real numbers $x_0,\cdots,x_k$ such that $0<x_{0}<x_{1}<\dots<x_{k}$ and
      $$u|_{|x|=x_j}\cdot u|_{|x|=x_{j+1}}<0,\,\,\,\,\,\,j=0,\dots,k-1.$$
Denote the nodal number of the function $u$ by $n(u)$.
\end{definition}

We always assume that the functions we discuss have finite nodal numbers.
\begin{definition}
For a continuous radial function $u$ with $n(u)=k$ and $u(x_{0})>0$, we define its $q$-th bump for $q=1,...,k+1$, by
      \begin{align}
      u_{1}(x)&=\chi_{\{u>0\}}\cdot\chi_{\{|x|<x_{1}\}}\cdot u(x),  \nonumber\\
      u_{q}(x)&=\chi_{\{(-1)^{q-1}u>0\}}\cdot\chi_{\{x_{q-2}<|x|<x_{q}\}}\cdot u(x),\,\,\,\,q=2,\dots,k,\nonumber\\
      u_{k+1}(x)&=\chi_{\{u(x_{k})\cdot u>0\}}\cdot\chi_{\{x_{q-1}<|x|\}}\cdot u(x).\nonumber
      \end{align}
      For a radial function $u$ with $n(u)=k$ and $u(x_{0})<0$, we define its $q$-th bump $q=1,...,k+1$ by
      \begin{align}
      u_{1}(x)&=\chi_{\{u<0\}}\cdot\chi_{\{|x|<x_{1}\}}\cdot u(x),  \nonumber\\
      u_{q}(x)&=\chi_{\{(-1)^{q-1}u<0\}}\cdot\chi_{\{x_{q-2}<|x|<x_{q}\}}\cdot u(x),\,\,\,\,q=2,\dots,k,\nonumber\\
      u_{k+1}(x)&=\chi_{\{u(x_{k})\cdot u>0\}}\cdot\chi_{\{x_{q-1}<|x|\}}\cdot u(x).\nonumber
      \end{align}

\end{definition}

\begin{remark}
To avoid confusion, for the $j$-th component $u_{j}$ of $U=(u_{1},\dots,u_{N})$, we denote its $q$-th bump by $u_{j,q}$.
\end{remark}

Now we present the concept of trivial solutions, semi-trivial solutions, and non-trivial solutions.
\begin{definition}
For a solution $U=(u_1,\cdots,u_N)$ to Problem (\ref{e:AAAA}), $U$ is said to be a trivial solution if for any $j=1,\cdots,N$, $u_j\equiv 0$ in $B_1$. $U$ is semi-trivial, if it is not trivial but there is a $j_0=1,\cdots,N$ with $u_{j_0}\equiv 0$ in $B_1$. $U$ is non-trivial if for any $j=1,\cdots,N$, $u_j$ is non-zero.
\end{definition}

\subsection{Liouville-type theorems and a priori estimates}
In this part, we aim to obtain a priori estimates for the solutions to Problem (\ref{e:AAAA}). Firstly, we prove Theorem \ref{t:Liouville}.

\noindent{\bf Proof of Theorem \ref{t:Liouville}.}
Our argument is analogues to the one in \cite[Proposition 21.2b]{QuitternBook2019}.
We argue by contradiction via a re-scaling method.
Let us denote
\begin{align}
\mathcal{U}_{\varepsilon,\beta}:=\big\{U\in(C(\mathbb{R}^3))^N\backslash\{0\}\big|U\mbox{ solves Problem }(\ref{e:CCCC})\mbox{ with }\varepsilon\leq\beta_{ij}\leq\beta\mbox{ for any }i\neq j\big\}.\nonumber
\end{align}
Here, $\beta>0$ is fixed in the assumptions. Assume that for a sequence of $\varepsilon_m\to0$ as $m\to0$, $\mathcal{U}_{\varepsilon_m,\beta}\neq\emptyset$. We can find a sequence $U_m\in \mathcal{U}_{\varepsilon_m,\beta}$. Denote $U_m=(u_{m,1},\cdots,u_{m,N})$ the solution to
\begin{equation}
    \left\{
   \begin{array}{lr}
     -{\Delta}u_j=\mu_j u_j^3+\sum_{i\neq j}\beta_{ij}^m u_i^2 u_j \,\,\,\,\,\,\,\,  \mbox{in }\mathbb{R}^3 ,\\
     u_j\in C^{2}(\mathbb{R}^3),\,u_j\mbox{ is radial,}\\
     u_j \mbox{is zero or changes its sign at most }P_j\mbox{ times and }(u_1,\cdots,u_N)\neq0,\\
     |u_j|_{L^\infty(\mathbb{R}^3)}\leq1\mbox{ for }j=1,\cdots,N
   \end{array}
   \right.
\end{equation}
with $\beta_{ij}^m\in[\varepsilon_m,\beta]$ for $i,j=1,\cdots,N$ and $i\neq j$.
Then, without loss of generality, we can assume that $u_{m,1}\neq 0$ and $\sup_{x\in\mathbb{R}^3}|u_{m,1}|=\sup_{j=1,\cdots,m}\sup_{x\in\mathbb{R}^3}|u_{m,j}|=:M_m$ holds for any $m=1,2,\cdots$. Up to a re-scaling, we can assume that $\frac{1}{2}\leq \sup_{x\in\mathbb{R}^3}|u_{m,1}|\leq 1$ for any $m$. Indeed, suppose $x_m\in\mathbb{R}^3$ such that $|u_{m,1}(x_m)|\geq\frac{M_m}{2}$.
Now, let us consider the sequence
\begin{align}
\widetilde{U}_m&:=\Big(M_m^{-1}U_{1,m}(M_m^{-1}(x+M_m x_m)),\cdots,M_m^{-1}U_{N,m}(M_m^{-1}(x+M_m x_m))\Big)\nonumber\\
&=:(\widetilde{U}_{m,1},\cdots, \widetilde{U}_{m,N}).\nonumber
\end{align}
satisfies $\frac{1}{2}\leq|\widetilde{U}_{m,1}(0)|\leq\sup_{j=1,\cdots,N}\sup_{x\in\mathbb{R}^3}|\widetilde{U}_{m,j}| =\sup_{x\in\mathbb{R}^3}|\widetilde{U}_{m,1}|\leq 1$.
If $\{M_m|x_m|\}_m$ is bounded,
as $m\to+\infty$, up to a translation in $\mathbb{R}^3$, there is a radial $C^{2,\alpha}$ vector-valued function $\widetilde{U}_\infty$ for some $\alpha\in(0,1)$ such that $\widetilde{U}_m\to\widetilde{U}_\infty$ in $C_{loc}^{2,\alpha}$. Writing $\widetilde{U}_\infty=\big(\widetilde{U}_{\infty,1},\cdots,\widetilde{U}_{\infty,N}\big)$, it follows that
\begin{align}\label{e:000}
-\Delta \widetilde{U}_{\infty,j}=\mu_j \widetilde{U}_{\infty,j}^3+\sum_{i\neq j}\beta_{ij}^\infty \widetilde{U}_{\infty,i}^2 \widetilde{U}_{\infty,j}\mbox{ in }\mathbb{R}^3
\end{align}
with $\lim_{m\to\infty}\beta_{ij}^m=\beta_{ij}^\infty\geq0$ for any $i,j=1,\cdots,N$ and $i\neq j$. $\widetilde{U}_{\infty,j}$ changes its sign at most $P_j$ times for $j=1,\cdots,N$. It should be noted that $\widetilde{U}_{\infty,1}\neq0$ and is bounded and nodes at most $P_1$ times. Then, denote $\mathcal{O}_1$ by the interior of support of the outmost bump of $\widetilde{U}_{\infty,1}$. Without loss of generality, we can assume that $\widetilde{U}_{\infty,1}|_{\mathcal{O}_1}>0$. Then,
\begin{align}
-\Delta \widetilde{U}_{\infty,1}\geq \mu_1 \widetilde{U}_{\infty,1}^3\mbox{ in }\mathcal{O}_1.
\end{align}
This contradicts with \cite[Theorem 1.3]{BidautVeron1989}.

If the sequence $\{M_m|x_m|\}_m$ is unbounded,
writing $r=|x-M_m x_m|$ and regarding $\widetilde{U}_{m,j}$ as one dimensional functions, then we have
\begin{align}
-\frac{d^2}{dr^2}\widetilde{U}_{m,j}-\frac{2}{r+M_m|x_{m}|}\frac{d}{dr}\widetilde{U}_{m,j}=\mu_j\widetilde{U}_{m,j}^3 +\sum_{i\neq j}\beta_{ij}^m\widetilde{U}_{m,i}^2 \widetilde{U}_{m,j}\mbox{ for }r\in(-M_m |x_m|,\infty)\nonumber
\end{align}
for any $j=1,\cdots,N$. Then, the $C_{loc}^{2,\alpha}$-limit of the sequence $\widetilde{U}_m$ satisfies
\begin{align}\label{e:intermiediate}
-\widetilde{U}_{\infty,j}''=\mu_j \widetilde{U}^3_{\infty,j}+\sum_{i\neq j} \beta_{ij}^\infty\widetilde{U}_{\infty,i}^2 \widetilde{U}_{\infty,j}\mbox{ for }r\in\mathbb{R}
\end{align}
for any $j=1,\cdots,N$. And $\widetilde{U}_{\infty,1}(0)\geq\frac{1}{2}$ and $\widetilde{U}_{\infty,1}$ has at most $P_1$ zeroes. Let us
denote $-\widetilde{U}_{\infty,1}'':=f(x,\widetilde{U}_{\infty,1})$. Using (\ref{e:intermiediate}), it is easy to verify that, $\widetilde{U}_{\infty,1}f(x,\widetilde{U}_{\infty,1})>0$ if $\widetilde{U}_{\infty,1}\neq 0$. By \cite[Theorem 1]{Wong1968}, $\widetilde{U}_{\infty,1}$ oscillates on $\mathbb{R}$.
We have a contradiction.

Therefore, the assumption is invalid and Theorem \ref{t:Liouville} holds.

\begin{flushright}
$\Box$
\end{flushright}

Two immediate consequences can be obtained.
\begin{corollary}\label{coro:nonnegativebeta}
For any $P_1\cdots,P_N\in\mathbb{N}$, if $\beta_{ij}\geq0$ for any $i,j=1,\cdots,N$ and $i\neq j$,
Problem (\ref{e:CCCC}) admits no non-trivial solution.
\end{corollary}

\begin{corollary}\label{coro:smallbeta}
For any $P_1\cdots,P_N\in\mathbb{N}$, there exists a positive number $\varepsilon_1=\varepsilon_1(\mu_1,\cdots,\mu_N;P_1,\cdots,P_N)>0$ such that if $|\beta_{ij}|<\varepsilon_1$, Problem (\ref{e:CCCC}) admits no non-trivial solution.
\end{corollary}

A one-dimensional analogue to the above results can be proved using a similar approach.
\begin{theorem}\label{t:Lioville2}
Consider the problem
\begin{equation}\label{e:DDDD}
    \left\{
   \begin{array}{lr}
     -u_j''=\mu_j u_j^3+\sum_{i\neq j}\beta_{ij} u_i^2 u_j \,\,\,\,\,\,\,\,  \mbox{in }\mathbb{R} ,\\
     u_j\in C^{2}(\mathbb{R}),\,u_j\mbox{ is radial,}\\
     u_j\mbox{ is zero or changes its sign at most }P_j\mbox{ times and }(u_1,\cdots,u_N)\neq0,\\
     |u_j|_{L^\infty(\mathbb{R})}\leq1\mbox{ for }j=1,\cdots,N.
   \end{array}
   \right.
\end{equation}
and the problem
\begin{equation}\label{e:EEEE}
    \left\{
   \begin{array}{lr}
     -u_j''=\mu_j u_j^3+\sum_{i\neq j}\beta_{ij} u_i^2 u_j \,\,\,\,\,\,\,\,  \mbox{in }[0,+\infty) ,\\
     u_j\in C^{2}([0,+\infty)),\,u_j\mbox{ is radial,}\\
     u_j\mbox{ is zero or changes its sign at most }P_j\mbox{ times and }(u_1,\cdots,u_N)\neq 0,\\
     |u_j|_{L^\infty([0,+\infty))}\leq1,\,u_j(0)=0\mbox{ for }j=1,\cdots,N.
   \end{array}
   \right.
\end{equation}
Here, $\mu_j>0$ are constants for $j=1,\cdots,N$.
For any $P_1\cdots,P_N\in\mathbb{N}$, and any $\beta\geq0$, there exists a positive number $\varepsilon_2=\varepsilon_2(\mu_1,\cdots,\mu_N;P_1,\cdots,P_N;\beta)>0$ such that if $\beta_{ij}\in[-\varepsilon_2,\beta]$ for any $i,j=1,\cdots,N$ and $i\neq j$,
Problem (\ref{e:DDDD}) and Problem (\ref{e:EEEE}) admit no non-trivial solution.
\end{theorem}

Similar to the above result, we also have two immediate corollaries.

\begin{corollary}\label{coro:nonnegativebeta1}
For any $P_1\cdots,P_N\in\mathbb{N}$, if $\beta_{ij}\geq0$ for any $i,j=1,\cdots,N$ and $i\neq j$,
Problem (\ref{e:DDDD}) and Problem (\ref{e:EEEE}) admit no non-trivial solution.
\end{corollary}

\begin{corollary}\label{coro:smallbeta1}
For any $P_1\cdots,P_N\in\mathbb{N}$, there exists a positive number $\varepsilon_3=\varepsilon_3(\mu_1,\cdots,\mu_N;P_1,\cdots,P_N)>0$ such that if $|\beta_{ij}|<\varepsilon_3$, Problem (\ref{e:DDDD}) and Problem (\ref{e:EEEE}) admit no non-trivial solution.
\end{corollary}

\begin{remark}\label{r:novelty}
The novelty of these Liouville-type theorems lies in the absence of the need for comparison between the components to conclude non-existence. However, comparisons are still crucial in a broader range of parameters. This is because the solutions to Problem (\ref{e:AAAA}) with a component-wise prescribed number of nodes may have no a priori estimate, especially for the case of $\beta_{ij}\equiv\beta$, $\mu_j\equiv\mu$ and $\beta<-\mu$ for any $i,j=1,\cdots,N$ and $i\neq j$, cf. \cite{LiWang2021,LiWangNEW,Quittner2021}. It is unreasonable to expect a Liouville-type theorem under these circumstance.
\end{remark}

\subsection{\cite[Theorem 1.3]{Tanaka2016} and the corresponding non-degeneracy result}
Consider the scalar field counterpart of Problem (\ref{e:AAAA}), i.e.,
\begin{equation}\label{e:main}
    \left\{
   \begin{array}{lr}
     -{\Delta}u+\lambda u=\mu u^3 \,\,\,\,\,\,\,\,  \mbox{in }B_\rho\subset\mathbb{R}^3 ,\\
     u\in H_{0,r}^1(B_\rho).
   \end{array}
   \right.
\end{equation}
Here, $\lambda,\mu>0$ are constant. $B$ is a ball of radius $\rho$ with the origin as its centre in $\mathbb{R}^3$. \cite[Theorem 1.3]{Tanaka2016} addresses that
\begin{theorem}
Assume that $\lambda=\mu=1$ and $\rho\leq 1$. For any $P=0,1,\cdots$, Problem (\ref{e:main}) admits an unique radial solution $u$ changing its sign exactly $P$ times with $u(0)>0$.
\end{theorem}
In this note, we would like to go slightly further on this direction, i.e., to prove the following result:
\begin{theorem}\label{t:main1}
Consider Problem (\ref{e:main}) with $\rho=1$. For any $P\in\mathbb{N}$, we have
\begin{itemize}
  \item [$(1).$] if $\lambda\geq1$, for any $k=0,1,\cdots$, Problem (\ref{e:main}) admits an unique radial solution $u$ changing its sign exactly $P$ times with $u(0)>0$;
  \item [$(2).$] the above solution is non-degenerate in the radial functional space $H_{0,r}^1(B_1)$.
\end{itemize}
\end{theorem}
\begin{remark}
Here, by the word \emph{non-degenerate in the radial functional space}, we mean that for any $\varphi\in H_{0,r}^1(B_1)$, we have
\begin{align}
-\Delta\varphi +\varphi =3 u^2\varphi.\nonumber
\end{align}
\end{remark}

Now we begin to prove Theorem \ref{t:main}.
~~

\noindent{\bf Proof of $(1)$.}
This computation follows a routine process as described in \cite{Tanaka2016}.
The only difference lies in the transformation of the radial functions defined on $B_1$ to one-dimensional functions, which we will restate as follows.
Consider
\begin{align}
-\frac{d^2}{dr^2}\phi+\frac{2}{r}\frac{d}{dr}\phi-\lambda \phi=0\nonumber
\end{align}
with
\begin{align}
\phi(r)=\frac{\sqrt{\lambda}\cosh(\frac{r}{\sqrt{\lambda}})}{r}.\nonumber
\end{align}
Then, we can define
\begin{align}
t=\int_0^r\frac{\lambda ds}{\big(\cosh(\frac{s}{\sqrt{\lambda}})\big)^2}\nonumber
\end{align}
and
\begin{align}
y(t)=\frac{u(r)}{\phi(r)}.\nonumber
\end{align}
It is easy to verify that
\begin{itemize}
  \item $\frac{d}{dt}y(t)=\frac{r^2}{\lambda}\cdot[\frac{d}{dr}u\phi - u\frac{d}{dr}\phi]$;
  \item $\frac{d^2}{dt^2}y+h(t)y=0$ with $h(t)=\lambda\mu^3\cdot\frac{(\cosh(\frac{r}{\sqrt{\lambda}}))^6}{r^2}$.
\end{itemize}
The rest part follows a routine similar to that of \cite[Corollary 3.4]{Tanaka2016}, which will not be discussed further.

~~

~~
\noindent{\bf Proof of $(2)$.}
The proof will be completed if we verify the following claim.
\begin{claim}
For the radial solution $u$ in Theorem \ref{t:main1}, problem
\begin{equation}\label{e:linearized}
    \left\{
   \begin{array}{lr}
     -{\Delta}v+\lambda v=3\mu u^2v \,\,\,\,\,\,\,\,  \mbox{in }B\subset\mathbb{R}^3 ,\\
     v\in H_{0,r}^1(B)
   \end{array}
   \right.
\end{equation}
admits only trivial solution, i.e., $v=0$.
\end{claim}
By the transform
\begin{align}
z(t)=\frac{v(r)}{\phi(r)},\nonumber
\end{align}
we have
\begin{equation}\label{e:linearized}
    \left\{
   \begin{array}{lr}
   \frac{d^2}{dt^2}z+3h(t)y^2z=0,\\
   z(0)=z(T)=0
   \end{array}
   \right.
\end{equation}
with $T=\int_0^1\frac{\lambda ds}{\big(\cosh(\frac{s}{\sqrt{\lambda}})\big)^2}$.
Here, since
\begin{align}
z(0)=\lim_{r\to0+}\frac{r\cdot v(r)}{\sqrt{\lambda}\cosh(\frac{r}{\sqrt{\lambda}})},\nonumber
\end{align}
we will have $z(0)=0$ as long as $v$ is a classical solution. $z(T)=0$ is due to $v\in H_{0,r}^1(B_1)$.
We intend to prove that $z\equiv 0$ in $[0,T]$. Arguing by contradiction.
Assume that $z\neq 0$. Since Problem (\ref{e:linearized}) is linear in $z$, we can assume $z'(0)=1$, without loss of generality. Consequently, the problem
\begin{equation}\label{e:linearized2}
    \left\{
   \begin{array}{lr}
   \frac{d^2}{dt^2}z+3h(t)y^2z=0,\\
   z(0)=0,\,z'(0)=1
   \end{array}
   \right.
\end{equation}
admits a non-trivial solution $z$ with $z(T)=0$. One the other hand, notice that
\begin{itemize}
  \item [$(1).$] $T$, the end point of the interval $[0,T]$, is also the largest zero of the function $z$;
  \item [$(2).$] Problem (\ref{e:linearized2}) admits an unique $C^2$-solution.
\end{itemize}
The first assertion is evident and the second is guaranteed by \cite[Proposition 4.1]{Tanaka2016}. Then, \cite[Lemma 5.2]{Tanaka2016} gives us that $z(T)\neq 0$. This is a impossible. Therefore,
$0$ is the only solution to Problem (\ref{e:linearized}) due to the its unique solvability.
This proves Theorem \ref{t:main1}.

\begin{flushright}
$\Box$
\end{flushright}

An immediate consequence of Theorem \ref{t:main1} concerns the invertibility of the second order derivative of the $C^2$ functional
\begin{align}\label{ENERGY}
I_{\vec{\beta}}(u_1,\cdots,u_N)=\frac{1}{2}\sum_{j=1}^N\int_{B_1}|\nabla u_j|^2 +\lambda_j u_j^2-\frac{1}{4}\sum_{j=1}^N\int_{B_1}\mu_j u_j^4+\sum_{i\neq j}\beta_{ij} u_i^2 u_j^2.
\end{align}
Here, $\vec{\beta}=(\beta_{ij})_{i\neq j;i,j=1,\cdots,N}$. Letting $u_{j,P_{j}}$ solve
$$-\Delta u+\lambda_j u=\mu_j u^3\mbox{ in }B_1$$
with $n(u_j)=P_j$. Then, we get
\begin{corollary}\label{coro:invertible}
It holds that the operator $I''|_{\vec{\beta}=0}(u_{1,P_1},\cdots,u_{N,P_N})$ is invertible in $(H_{0,r}^1(B_1))^N$.
\end{corollary}
This is a direct consequence of the $C^2$-continuity of the functional $I_{\vec{\beta}}$.

\section{Proof of Theorem \ref{t:main} and Theorem \ref{t:boundedness}}
This section presents the proof for Theorem \ref{t:main} and Theorem \ref{t:boundedness}. The former theorem is an immediate consequence of the a priori estimate (see Theorem \ref{t:boundedness}), the uniqueness (see Theorem \ref{t:main1}) and non-degeneracy result (see Corollary \ref{coro:invertible}) through an argument in \cite{Ikoma2009}. To this end, we need to verify Theorem \ref{t:boundedness} in advance.

\subsection{Proof of Theorem \ref{t:boundedness} and its corollaries}
The proof of Theorem \ref{t:boundedness} is a routine and we only outline it. Arguing by contradiction, suppose for a $\beta>0$, we can find a sequence of numbers $\beta_{ij}^{(m)}\in[-\varepsilon_*,\beta]$ for any $i,j=1,\cdots,N$ and $i\neq j$ with
\begin{itemize}
  \item $\varepsilon_*=\min(\varepsilon_0,\varepsilon_2)$, the constants in Theorem \ref{t:Liouville} and Theorem \ref{t:Lioville2};
  \item a sequence of functions $(u_1^{(m)},\cdots,u_N^{(m)})$ with
    \begin{equation}
    \left\{
   \begin{array}{lr}
     -{\Delta}u_j^{(m)}+\lambda_j u_j^{(m)}=\mu_j (u_j^{(m)})^3+\sum_{i\neq j}\beta^{(m)}_{ij} (u_i^{(m)})^2 u_j^{(m)} \,\,\,\,\,\,\,\,  \mbox{in }B_1 ,\nonumber\\
     u_j^{(m)}\in H_{0,r}^1(B_1),\,n(u_j^{(m)})=P_j\mbox{ for }j=1,\cdots,N;\nonumber
   \end{array}
   \right.
\end{equation}
  \item without loss of generality, $|u_1^{(m)}|_\infty=\max_{j=1,\cdots,N}|u_j^{(m)}|_\infty\to\infty$ as $m\to\infty$.
\end{itemize}
As in Section 2.2, denote $M'_m=|u_1^{(m)}|_\infty$ and $x'_m\in B_1$ such that $|u_1^{(m)}(x'_m)|=|u_1^m|_\infty$ for any $m\in\mathbb{N}$. Consider the re-scaling functions
\begin{align}
\hat{U}_m&:=\Big((M'_m)^{-1}u_1^{(m)}((M'_m)^{-1}(x+M'_m x'_m)),\cdots,(M'_m)^{-1}u_N^{(m)}((M'_m)^{-1}(x+M'_m x'_m))\Big)\nonumber\\
&=:(\hat{U}_{m,1},\cdots, \hat{U}_{m,N}).\nonumber
\end{align}
If the sequence $\{M'_m|x'_m|\}_m$ is bounded, the $C^{2,\alpha}_{loc}$ limit of $\hat{U}_m$ will be a non-trivial classical solution to Problem (\ref{e:CCCC}), which contradicts Theorem \ref{t:Liouville}. If the sequence $\{M'_m|x'_m|\}_m$ is unbounded, the limit will tend to a non-trivial classical solution to Problem (\ref{e:DDDD}) or Problem (\ref{e:EEEE}). This contradicts Theorem \ref{t:Lioville2}. Therefore, Theorem \ref{t:boundedness} is proved.
\begin{flushright}
$\Box$
\end{flushright}

As we did in Section 2.2, two corollaries follow immediately,
\begin{corollary}\label{coro:nonnegativebeta2}
For any $P_1\cdots,P_N\in\mathbb{N}$, if $\beta_{ij}\geq0$ for any $i,j=1,\cdots,N$ and $i\neq j$, there exists a constant $C=C(\lambda_1,\cdots,\lambda_N;\mu_1,\cdots,\mu_N;\beta_{ij};P_1,\cdots,P_N)>0$ such that for any solution $U$ to
Problem (\ref{e:AAAA}), we have $|U|_\infty\leq C$.
\end{corollary}

\begin{corollary}\label{coro:smallbeta2}
For any $P_1\cdots,P_N\in\mathbb{N}$, there exist positive numbers $\varepsilon_4=\varepsilon_4(\mu_1,\cdots,\mu_N;P_1,\cdots,P_N)>0$ and $C=C(\lambda_1,\cdots,\lambda_N;\mu_1,\cdots,\mu_N;P_1,\cdots,P_N)>0$ such that if $|\beta_{ij}|<\varepsilon_4$, for any solution $U$ to Problem (\ref{e:AAAA}), we have $|U|_\infty\leq C$.
\end{corollary}

\subsection{Proof of Theorem \ref{t:main}}
We argue by contradiction.
Suppose that there is a sequence $(\beta^m_{ij})_{i\neq j,i,j=1,\cdots,N;m=1,2,\cdots}$ such that $\beta_{ij}^m\to0$ as $m\to+\infty$ for any $i\neq j$ and $i,j=1,\cdots,N$ and there are two different solutions $U_m^1$ and $U_m^2$ to the problem
  \begin{equation}
    \left\{
   \begin{array}{lr}
     -{\Delta}u_j+\lambda_j u_j=\mu_j u_j^3+\sum_{i\neq j}\beta^m_{ij} u_i^2 u_j \,\,\,\,\,\,\,\,  \mbox{in }B_1 ,\nonumber\\
     u_j\in H_{0,r}^1(B_1),\,n(u_j)=P_j\mbox{ for }j=1,\cdots,N.\nonumber
   \end{array}
   \right.
\end{equation}
Here, $\lambda_j\geq1$ and $\mu_j>0$ for $j=1,\cdots,N$. By Corollary \ref{coro:smallbeta2}, there exists a positive constant $C$ independence in $m$ such that $|U_1^m|_{L^\infty(B_1)},|U_2^m|_{L^\infty(B_1)}\leq C$.
Then, for any $p\geq2$, a $W^{2,p}$-estimate yields that $\|U_1^m\|_{W^{2,p}(B_1)},\|U_2^m\|_{W^{2,p}(B_1)}\leq C$. Then, up to a subsequence, there are $U_{1}$ and $U_{2}$ such that
\begin{align}
U_1^m\to U_1\mbox{ and }U_2^m\to U_2\mbox{ in }(H_0^1(B_1))^N\nonumber
\end{align}
as $m\to \infty$. And it is evident that $U_1$ and $U_2$ solve
\begin{equation}
    \left\{
   \begin{array}{lr}
     -{\Delta}u_j+\lambda_j u_j=\mu_j u_j^3 \,\,\,\,\,\,\,\,  \mbox{in }B_1 ,\nonumber\\
     u_j\in H_{0,r}^1(B_1),\,n(u_j)\leq P_j\mbox{ for }j=1,\cdots,N.\nonumber
   \end{array}
   \right.
\end{equation}
Now we will verify that for any $j=1,\cdots,N$ we have $n(u_j)= P_j$. Denote $U_i^m=(u_{i,1}^m,\cdots,u_{i,N}^m)$ for $i=1,2$. Therefore, $n(u_{i,j}^m)=P_j$. Denote by $w^m$ a bump of $u_{i,j}^m$. Multiplying the $j$-th equation of Problem (\ref{e:AAAA}) and integrating over $B_1$, we get
\begin{align}
\int_{B_1}|\nabla w^m|^2+\lambda_j |w^m|^2 & =\mu_j\int_{B_1}|w^m|^4+\sum_{i\neq j}\beta_{ij}\int_{B_1} u_i^2 |w^m|^2\nonumber\\
&\leq\mu_j|w^m|_4^4+C_0\cdot\sup_{i\neq j}|\beta_{ij}|\cdot|w^m|_4^2.\nonumber
\end{align}
Using Sobolev's inequality,
\begin{align}
C_1|w^m|_4^2\leq \mu_j|w^m|_4^4+C\cdot\sup_{i\neq j}|\beta_{ij}|\cdot|w^m|_4^2.\nonumber
\end{align}
Letting $\sup_{i\neq j}|\beta_{ij}|<\frac{C_1}{2C_0}$, $|w^m|_4\geq C$ for some constant $C>0$ independence in $m$. The $L^4$ convergence yields the conclusion that $n(u_j)= P_j$. By Theorem \ref{t:main1}, we have
$U_1=U_2$, i.e., $\|U_1^m - U_2^m\|\to 0$ as $m\to\infty$.

On the other hand, due to the implicit function theorem (cf. \cite[Theorem 1.2.1]{Chang2005BOOK}), there is a constant $b>0$ such that for any set of numbers $\vec{\beta}=(\beta_{ij})_{i,j=1,\cdots,N;i\neq j}$ with $|\beta_{ij}|<b$, the equation
\begin{align}
I'_{\vec{\beta}}(u_1,\cdots,u_N)=0\nonumber
\end{align}
admits an unique non-trivial solution in $\{U=(u_1,\cdots,u_N)|u_j\in H_{0,r}^1(B_1)\mbox{ and }\|u_j - u_{j,P_j}\|\leq \delta\}$ for some small $\delta>0$. This contradicts the convergence of $U_1^m$ and $U_2^m$. The assumption fails and Theorem \ref{t:main} is proved. The non-degeneracy mentioned in Remark \ref{r:nondegeneracy} is an immediate consequence of the $C^2$ continuity.

\begin{flushright}
$\Box$
\end{flushright}

%%%%%%%%%%%%%%%%%%%%%%%%%%%%%%%%%%%%%%%%%%%%%%%%%%%%%%%%%%%%%%%%%%%%%%%%%%%%%%%%%%%%%%%%%%%%%%%%%%%%%%%%%%%%%%%%%%%%%%%%%%%%%%%%%%%%%%%%%%%%
%%%%%%%%%%%%%%%%%%%%%%%%%%%%%%%%%%%%%%%%%%%%%%%%%%%%%%%%%%%%%%%%%%%%%%%%%%%%%%%%%%%%%%%%%%%%%%%%%%%%%%%%%%%%%%%%%%%%%%%%%%%%%%%%%%%%%%%%%%%%

\section*{Acknowledgements}
\addcontentsline{toc}{section}{Acknowledgements}
Both of the authors would like to express their sincere appreciation to the anonymous reviewers for their helpful suggestions on improving this paper. Their thank also goes to Professor Pavol Quittner, Professor Satoshi Tanaka and Professor Zhi-Qiang Wang for their various helpful communications.
In this research, HL was supported by FAPESP PROC 2022/15812-0 and OHM was supported by CNPq PROC 303256/2022-2 AND FAPESP PROC 2022/16407-1.

%%%%%%%%%%%%%%%%%%%%%%%%%%%%%%%%%%%%%%%%%%%%%%%%%%%%%%%%%%%%%%%%%%%%%%%%%%%%%%%%%%%%%%%%%%%%%%%%%%%%%%%%%%%%%%%%%%%%%%%%%%%%%%%%%%%%%%%%%%%%
%%%%%%%%%%%%%%%%%%%%%%%%%%%%%%%%%%%%%%%%%%%%%%%%%%%%%%%%%%%%%%%%%%%%%%%%%%%%%%%%%%%%%%%%%%%%%%%%%%%%%%%%%%%%%%%%%%%%%%%%%%%%%%%%%%%%%%%%%%%%

\vspace{-0.11cm}

{

}


{\footnotesize

\begin {thebibliography}{44}

\bibitem{AoWeiYao2015}
Ao, W, Wei, J, Yao, W,
Uniqueness and nondegeneracy of sign-changing radial solutions to an almost critical elliptic problem.
October, 2015. Advances in Differential Equations 21(12).
%https://doi.org/10.4208/ata.OA-0009







\bibitem{BidautVeron1989}
Bidaut-V\'eron, M.-F,
Local and global behavior of solutions of quasilinear equations of Emden-Fowler type.
Arch. Rational Mech. Anal. 107(1989), no.4, 293-324.


\bibitem{Chang2005BOOK}
Chang, K.-C,
Methods in nonlinear analysis.
Springer Monographs in Mathematics. Berlin: Springer. ix, 439 p. (2005).



\bibitem{DancerWeiWeth2010}
Dancer, E. N, Wei, J, Weth, T,
A priori bounds versus multiple existence of positive solutions for a nonlinear Schr\"odinger system.
Annales de l'Institut Henri Poincare (C) Non Linear Analysis, 2010, 27(3) pp. 953-969.
%\url{https://doi.org/10.1016/j.anihpc.2010.01.009}

\bibitem{GidasNiNirenberg1979}
Gidas, B, Ni, W. M, Nirenberg, L,
Symmetry and related properties via the maximum principle.
Communications in Mathematical Physics, 1979, 68(3):209-243.
%\url{https://doi.org/10.1007/BF01221125}






\bibitem{GidasSpruck1981}
Gidas, B, Spruck, J,
A priori bounds for positive solutions of nonlinear elliptic equations.
Comm. Partial Differential Equations 6 (1981), no. 8, 883-901.
%\url{https://doi.org/10.1080/03605308108820196}




\bibitem{Ikoma2009}
Ikoma, N,
Uniqueness of positive solutions for a nonlinear elliptic system.
NoDEA Nonlinear Differential Equations Appl. 16 (2009), no. 5, 555-567.
%\url{https://doi.org/10.1007/s00030-009-0017-x}


\bibitem{KabeyaTanaka1999}
Kabeya, Y, Tanaka, K,
Uniqueness of positive radial solutions of semilinear elliptic equations in $\mathbb{R}^N$ and S\'er\'e's non-degeneracy condition.
Comm. Partial Differential Equations24(1999), no.3-4, 563-598.
%\url{https://doi.org/10.1080/03605309908821434}


\bibitem{KormanLiSchmidt2012}
Korman, P, Li, Y, Schmidt, D. S,
A computer assisted study of uniqueness of ground state solutions.
Journal of Computational $\And$ Applied Mathematics, 2012, 236(11):2838-2843.
%\url{https://doi.org/10.1016/j.cam.2012.01.020}


\bibitem{Kwong1989}
Kwong, M. K,
Uniqueness of positive solutions of $\Delta u - u + u^p =0$ in $\mathbb{R}^n$.
Archive for Rational Mechanics and Analysis, 1989, 105(3), 243-266.
%\url{https://doi.org/10.1007/BF00251502}


\bibitem{LiWang2021}
Li, H, Wang, Z.-Q,
Multiple nodal solutions having shared componentwise nodal numbers for coupled Schr\"odinger equations.
Journal of Functional Analysis 280(7), 2021:108872.
%\url{https://doi.org/10.1016/j.jfa.2020.108872.}



%\bibitem{LiWangNEW}
%Li, H, Wang, Z.-Q,
%Multiple nodal solutions having shared componentwise nodal numbers for coupled Schr\"odinger equations. II.
%In preparation.


\bibitem{LiuLiuChang2015}
Liu, H, Liu, Z, Chang, J,
Existence and uniqueness of positive solutions of nonlinear Schr\"odinger systems.
Proceedings of the Royal Society of Edinburgh: Section A Mathematics, 145(2), 2015:365-390.
%\url{https://doi.org/10.1017/S0308210513000711}


\bibitem{LiuWang2019}
Liu, Z, Wang, Z.-Q,
Vector Solutions with Prescribed Component-Wise Nodes for a Schr\"odinger System.
Analysis in Theory and Applications 35(3), 2019:288-311.
%\url{https://doi.org/10.4208/ata.OA-0009}









\bibitem{Quittner2021}
Quittner, Pavol
Liouville theorem and a priori estimates of radial solutions for a non-cooperative elliptic system.
Nonlinear Anal., Theory Methods Appl., Ser. A, Theory Methods. 222, 2022:112971.
%\url{https://doi.org/10.1016/j.na.2022.112971}




\bibitem{QuitterSouplet2012}
Quittner, P, Souplet, P,
Optimal Liouville-type Theorems for Noncooperative Elliptic Schr\"odinger Systems and Applications.
Communications in Mathematical Physics. 311, 2012:1-19.
%\url{https://doi.org/10.1007/s00220-012-1440-0}


\bibitem{QuitternBook2019}
Quittner, P, Souplet, P,
Superlinear Parabolic Problems: Blow-up, Global Existence and Steady States. Second edition.
Birkh\"auser Basel. 2019.



\bibitem{QuittnerCommunication}
Quittner, P,
Personal communications. 2023.

\bibitem{Tanaka2016}
Tanaka, S,
Uniqueness of sign-changing radial solutions for $\Delta u-u +|u|^{p-1}u=0$ in some ball and annulus.
Journal of Mathematical Analysis and Applications, 439(1), 2016:154-170.
%\url{https://doi.org/10.1016/j.jmaa.2016.02.036}


\bibitem{YaoWei2011}
Wei, J, Yao, W,
Uniqueness of positive solutions to some coupled nonlinear Schr\"odinger equations.
Communications on Pure $\And$ Applied Analysis, 2011, 11(3):1003-1011.
%\url{https://doi.org/10.3934/cpaa.2012.11.1003}


\bibitem{Wong1968}
Wong, J. S. W,
On second order nonlinear oscillation.
Funkcial. Ekvac.11(1968), 207-234.


\bibitem{ZhouWang2020}
Zhou, L, Wang, Z.-Q,
Uniqueness of positive solutions to some Schr\"odinger systems.
Nonlinear Analysis, 2020, 195(2):111750.
%\url{https://doi.org/10.1016/j.na.2020.111750}



\end {thebibliography}
}
\end{document}